\documentclass[letterpaper,11pt,reqno]{amsart}
\usepackage[utf8]{inputenc}

\usepackage[
	theoremdefs,
	final
]{latexdev}
\usepackage[numbers]{natbib} 
\usepackage[all]{xy}
\usepackage{booktabs}
\usepackage{amssymb}
\usepackage{color}
\usepackage{amsxtra} 
\usepackage{tikz}

\usepackage{mathrsfs} 

\newcommand{\arxivorpnas}[2]{#1}

\usepackage{amsmath,amssymb}

\DeclareMathOperator{\divv}{div}

\newcommand{\ee}{\mathrm{e}}
\newcommand{\ii}{\mathrm{i}}
\newcommand{\pair}[1]{\left\langle #1 \right\rangle}

\newcommand{\inner}[1]{\langle\!\langle #1 \rangle\!\rangle}
\providecommand{\norm}[1]{\lVert#1\rVert}
\providecommand{\abs}[1]{\lvert#1\rvert}

\newcommand{\ud}{\mathrm{d}}

\newcommand{\RR}{{\mathbb R}}
\newcommand{\CC}{{\mathbb C}}

\newcommand{\vol}{\mu}
\newcommand{\met}{\mathsf{g}}
\newcommand{\Diff}{\mathrm{Diff}}
\newcommand{\Xcal}{\mathfrak{X}}
\newcommand{\Diffvol}{{\Diff_\vol}}

\newcommand{\LieD}{\mathcal{L}}

\newcommand{\Dens}{\mathrm{Dens}}

\newcommand*\id{\mathrm{id}}

\newcommand{\Met}{\mathsf{G}}
\newcommand{\MetW}{ {\bar{\mathsf{G}}} }
\newcommand{\MetF}{ {\bar{\mathsf{G}}} }

\newcommand{\sslash}{\mathbin{\mkern-3mu/\mkern-5mu/\mkern-2mu}}

\newcommand{\marginnote}[1]
{
}

\newcounter{gm}

\newcounter{bk}

\newcounter{km}

\title[Geometric Hydrodynamics via Madelung Transform]
{Geometric Hydrodynamics via Madelung Transform}\author{Boris Khesin \and Gerard \texorpdfstring{Misio\l ek}{Misiolek} \and Klas Modin}
\date{\today}

\begin{document} 
\begin{abstract} 
We introduce a geometric framework to study Newton's equations on infinite-dimensional configuration 
spaces of diffeomorphisms and smooth probability densities. It turns out that several important PDEs 
of hydrodynamical origin can be described in this framework in a natural way. 
In particular, the Madelung transform between the Schr\"odinger equation and 
Newton's equations is a symplectomorphism of the corresponding phase spaces. 
Furthermore, the Madelung transform turns out to be a K\"ahler map when the space of densities 
is equipped with the Fisher-Rao information metric. 
We describe several dynamical applications of these results.
\end{abstract} 

\maketitle 

\tableofcontents 


\newcommand{\eqtable}{
\begin{table}
	\arxivorpnas{}{
	\caption{Examples of Newton's equations.}\label{tab:equations}
	}
	\centering
	\begin{tabular}{ll}
		\toprule
		Wasserstein-Otto geometry & Fisher-Rao geometry \\
		\midrule
		\multicolumn{2}{c}{Newton's equations on $\Diff(M)$} \\[0.5ex]
		$\bullet$ Classical mechanics & $\bullet$ $\mu$-Camassa-Holm \\
		$\bullet$ Burgers' inviscid & $\bullet$ Optimal information transport \\
		$\bullet$ Barotropic inviscid fluid & \\[1ex]
		\multicolumn{2}{c}{Newton's equations on $\Dens(M)$ or $T^*\Dens(M)$} \\[0.5ex]
		$\bullet$ Hamilton-Jacobi & $\bullet$ $\infty$-dim Neumann problem \\
		$\bullet$ Linear Schrödinger & $\bullet$ Klein-Gordon \\
		$\bullet$ Non-linear Schrödinger & $\bullet$ 2-component Hunter-Saxton \\
		$\bullet$ Vortex filament equation &  \\
		\bottomrule
	\end{tabular}
	\arxivorpnas{
	\caption{Examples of Newton's equations.}\label{tab:equations}
	}{}
\end{table}
}

\arxivorpnas{\eqtable}{}

\arxivorpnas{
\section{Introduction}\label{sec:intro}
}{}

\arxivorpnas{In}{\dropcap{I}n} a seminal 1966 paper, Arnold~\cite{Ar1966} showed that the Euler equations of an inviscid incompressible fluid 
can be reformulated as geodesic equations on an infinite-dimensional manifold of diffeomorphisms of a fluid domain. 
That paper was a foundation stone for a new branch of mathematics called \emph{geometric and topological hydrodynamics} 
\cite{ArKh1998} 
and many important PDEs of mathematical physics have been shown to fit Arnold's framework since. 
Examples include the Korteweg-de~Vries equation, the Camassa--Holm equation, 
magneto-hydrodynamics, \arxivorpnas{the Hunter--Saxton equation, and the Heisenberg spin chain.}{and the Hunter--Saxton equation.} 

Arnold's reformulation of the Euler equations in an elegant differential-geometric language allowed an insight into 
both analysis and geometry of the equations of fluid dynamics. 
For example, the sectional curvature of the group of diffeomorphisms of the fluid domain influences the dynamics of 
fluid motions via the equations of geodesic deviation, which had applications to hydrodynamic stability 
\cite{Ar1966, ArKh1998}.
Furthermore, a detailed study of the analytic properties of the associated Riemannian 
exponential map, begun by \citet{EbMa1970}, led to sharp local well-posedness results for the Cauchy problem 
of the Euler equations. 
One may expect that further study will shed new light on challenging problems of fluid dynamics, such
as regularity and persistence of solutions in 3D flows or the problem of fluid turbulence. 

In this paper we propose an extension of this approach to the case of \textit{Newton's equations} as a natural next step 
in Arnold's program. 
We are interested in second order equations that formally can be written as 
\begin{equation}\label{eq:newton_general} 
	\nabla_{\dot q}\dot q = - \nabla U(q), 
\end{equation} 
where $\nabla$ is the covariant derivative with respect to a certain Riemannian metric and $U$ is a potential function, cf. \cite{Sm1979}.
We develop a geometric framework for \arxivorpnas{the equations}{} \eqref{eq:newton_general} on the space of diffeomorphisms of a compact manifold.
Using infinite-dimensional Riemannian submersion techniques we show that these equations are closely related to 
Newton's equations on the space of smooth probability densities of the same underlying manifold. 
In particular, this is the case for the equations of compressible fluids. These equations have been long known to have a Hamiltonian formulation on the dual of 
a semi-direct product Lie algebra \cite{MaRaWe1984b}, while 
their Lagrangian Arnold-type formulation was best understood in terms of Newton's equations \cite{Sm1979}, as
its  Lagrangian in the semidirect product was not quadratic and led to additional difficulties \cite{HoMaRa1998}.
It turns out, however, that in  our framework of Riemannian submersions these equations have the following simple description\arxivorpnas{ (Section \ref{sub:compressible_euler_equations})}{}.
\begin{theorem}\label{thm:potential_compress}
The equation for potential solutions of compressible fluid in a compact domain is Newton's equation 
on the space of smooth probability densities with a potential function given by the fluid's internal energy. 
\end{theorem}
\arxivorpnas{}{\eqtable}
The proposed framework also reveals some unexpected connections between various results 
in fluid dynamics, optimal transport, information geometry 
and equations of mathematical physics, 
such as the Schrödinger equation, the Klein-Gordon equation, the Hunter-Saxton equation and its variants (see \autoref{tab:equations}).  
For instance, the classical Laplace eigenproblem 
can be seen as the problem of determining stationary solutions of a Fisher-Rao-Newton equation
 on the space of densities, 
which in turn describes geodesics on an infinite-dimensional ellipsoid through its formulation 
as a Neumann problem\arxivorpnas{ (Section \ref{sect:Neumann})}{}.

An important tool in our constructions is the Madelung transform which turns out to have a number of surprising properties and can be viewed as a symplectomorphism, an isometry, a K\"ahler map or a generalization of the Hasimoto transform depending on the context\arxivorpnas{ (\autoref{sec:madelung})}{}. 
Our study reveals that the geometric features of the Madelung transform 
are best understood not in the setting of the $L^2$-Wasserstein geometry 
but the Fisher-Rao geometry---the canonical Riemannian structure on the space of probability densities\arxivorpnas{ (Section \ref{sub:kahler_properties_of_madelung})}{}.
\begin{theorem}\label{thm:madelung_intro}
	The Madelung transform is a K\"ahler morphism between the cotangent bundle of the space of smooth probability densities, equipped with the (Sasaki)-Fisher-Rao metric, and an open subset (in the Fr\'echet topology of smooth functions) of the complex projective space of smooth wave functions, equipped with the Fubini-Study (or Bures) metric.
\end{theorem}
This result uncovers some surprising new links between hydrodynamics, quantum information geometry 
and geometric quantum mechanics. 

For conceptual clarity and brevity of the exposition we focus here on the formal aspects of 
the various infinite-dimensional geometric constructions.
Proofs of the theorems and a suitable functional-analytic setting of Sobolev and Fr\'echet spaces will be described in our forthcoming paper.

\arxivorpnas{
\bigskip

{\bf Acknowledgements:} A part of this work was completed while B.K.\ held
a Weston Visiting Professorship at the Weizmann Institute of Science. B.K.\ is grateful to the Weizmann
Institute for its support and kind hospitality. B.K.\ was also partially supported by an
NSERC research grant and a Simons fellowship.
K.M.\ was supported by EU Horizon 2020 grant No 661482, and by the Swedish Foundation for Strategic Research grant ICA12-0052. 
A part of this work was done while G.M. held the Ulam Chair visiting Professorship in University of Colorado at Boulder.
}{}
 
\arxivorpnas{
\section{Wasserstein geometry of the space of densities} \label{sec:wasserstein}
}{
\section*{Wasserstein geometry of the space of densities}
}

In this section we recall the main notions of the Wasserstein geometry on the space of diffeomorphisms $\Diff(M)$ 
and the space of smooth probability densities $\Dens(M)$\arxivorpnas{ (viewed as infinite-dimensional Fr\'echet manifolds)}{}
and we describe Newton's equations on these spaces. 

Let $(M,\met)$ be a compact Riemannian manifold with volume form $\vol$. 
Define an $L^2$-metric on $\Diff(M)$ by 
\begin{equation}\label{eq:L2met} 
\Met_{\varphi}(\dot\varphi,\dot\varphi) = \int_M \underbrace{\abs{\dot\varphi}^2}_{\met_{\varphi}(\dot\varphi,\dot\varphi)} \vol.
\end{equation}
Given a $C^1$ function $U\colon\Diff(M)\to\RR$ (a potential), Newton's equations on $\Diff(M)$ 
can be formally written as in \eqref{eq:newton_general}. 
%
We are interested in the special case in which the potential functions are of the form
%
\begin{equation}\label{eq:Ubar}
	U(\varphi) = \bar U(\varphi_*\vol),
\end{equation}
where $\bar U\colon \Omega^n(M)\to\RR$ and $\varphi_*\vol$ denotes the pushforward of the volume form $\mu$ 
by the diffeomorphism $\varphi$ in $\Diff(M)$.
Newton's equation on $\Diff(M)$ then takes the following form.
\begin{theorem}[cf.\ \cite{Sm1979}]\label{thm:newton_for_Ubar}
	Newton's equation on $\Diff(M)$ for the $L^2$-metric~\eqref{eq:L2met} and with a potential function \eqref{eq:Ubar} 
	has the form
	\begin{equation}\label{eq:newton_for_ubar}
		\nabla_{\dot\varphi}\dot\varphi = -\nabla\left( \frac{\delta \bar U}{\delta \varrho}(\varphi_*\vol) \right)\circ\varphi,
	\end{equation}
	where $ \varrho=\rho\mu$. 
	In the reduced variables $u=\dot\varphi\circ\varphi^{-1}$ and $\rho = \mathrm{det}(D\varphi^{-1})$ 
	equation \eqref{eq:newton_for_ubar} assumes the form
\begin{equation}\label{eq:reduced_newton_diff}
\left\{
  \begin{array}{l}
		~\dot u + \nabla_u u + \nabla  \dfrac{\delta \bar U}{\delta\varrho}(\varrho) = 0 \\\
		\dot\rho + \divv(\rho u) = 0.
\end{array} \right. 
\end{equation}
\end{theorem}
\arxivorpnas{Equations }{}\eqref{eq:reduced_newton_diff} admit an invariant subset of potential solutions $u=\nabla\theta$,
where $\theta\in C^\infty(M)$.
We now describe the geometric origin of this observation. 

\arxivorpnas{
\subsection{Riemannian submersion to densities} 
\label{sub:fibration}
}{
\subsection*{Riemannian submersion to densities}
}

The space of smooth probability densities on $M$ is the space of volume forms 
with total volume 1, namely 
\begin{equation*}
	\Dens(M) = \Big\{ \varrho\in \Omega^n(M)\mid \varrho > 0, \; \int_M \varrho = 1 \Big\}.
\end{equation*}
Since $\Dens(M)$ is an open subset of codimension one in an affine subspace of $\Omega^n(M)$,
the tangent bundle of $\Dens(M)$ is trivial:
$T\Dens(M)=\Dens(M)\times\Omega^n_0(M)$, 
where 
$\Omega^n_0(M) = \{ \alpha\in\Omega^n(M) \mid \int_M \alpha = 0 \}$.
Likewise, the (smooth part of the) cotangent bundle $T^*\Dens(M)$ is $\Dens(M)\times C^\infty(M)/\RR$.

Alternatively, $\Dens(M)$ can be viewed as the space of left cosets of the subgroup $\Diffvol(M)$ of 
volume-preserving diffeomorphisms of $M$ with the push-forward map $\pi (\varphi) = \varphi_\ast \vol$ 
defining a natural (left coset) projection $\pi : \Diff(M) \to \Dens(M)$. 
To take full advantage of this setup it is useful to equip the base space with a Sobolev $H^{-1}$-type metric 
which arises in the context of optimal mass-transport problems, cf.\ \cite{BeBr2000}. 

\begin{figure}
	\centering
	%
	\begin{tikzpicture}
		\node[anchor=south west, inner sep=0] (image) at (0,0) {\includegraphics[width=\arxivorpnas{0.5\textwidth}{0.32\textwidth}]{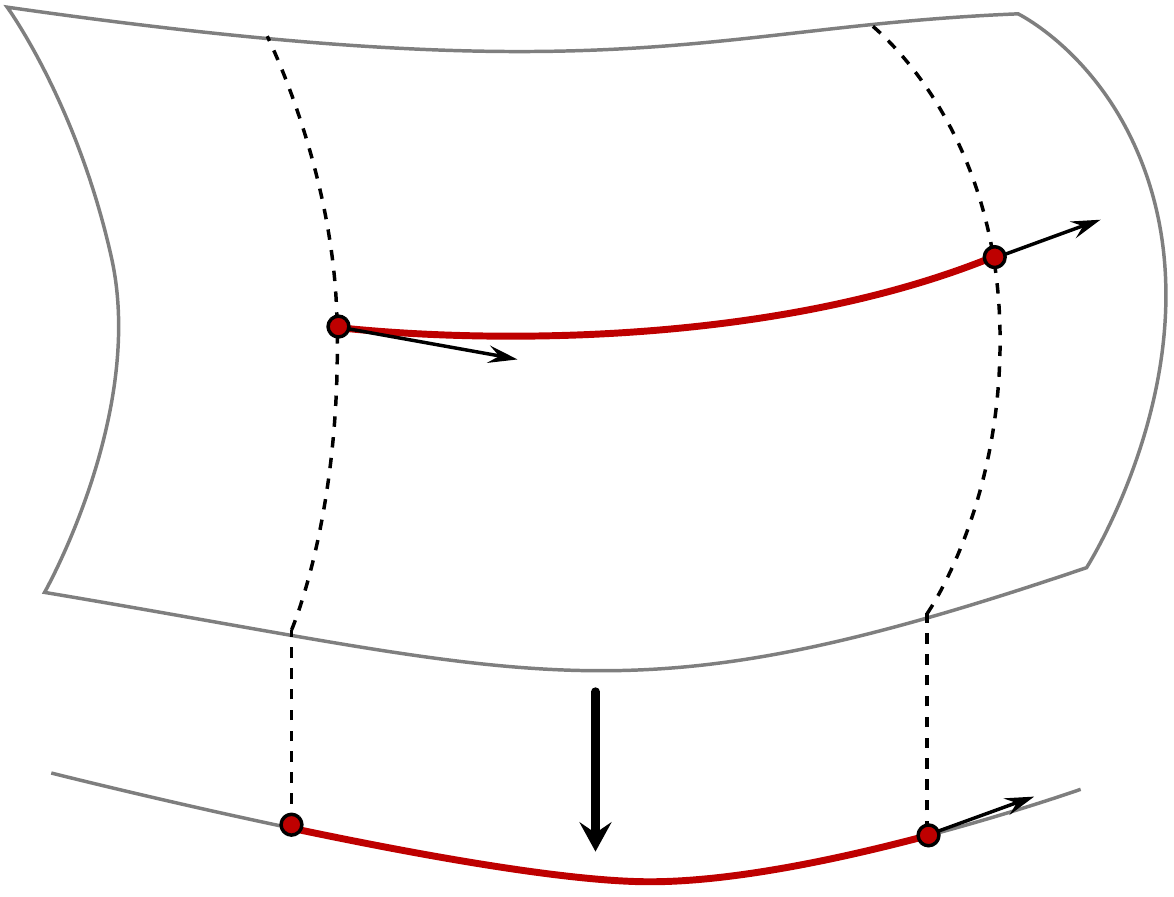}};
		\begin{scope}[x={(image.south east)},y={(image.north west)}]
			\coordinate (left_fiber) at (0.29,0.29) {};
			\coordinate (left_diffmu) at (0.21,0.3) {};
			\coordinate (right_fiber) at (0.84,0.33) {};
			\coordinate (hor_geo) at (0.57,0.64) {};
			\coordinate (geo) at (0.45,0.04) {};
			\coordinate (projection) at (0.51,0.16) {};
			\coordinate (diff) at (0.1,0.95) {};
			\coordinate (dens) at (0.1,0.1) {};
			\coordinate (id) at (0.29,0.64) {};
			\node[left, rotate=0] at (id) {$\mathrm{id}$};
			\coordinate (mu) at (0.245,0.07) {};
			\node[below, rotate=0] at (mu) {$\mu$};
			\coordinate (rho) at (0.79,0.065) {};
			\node[below, rotate=0] at (rho) {$\varrho$};
			\coordinate (dotrho) at (0.86,0.155) {};
			\node[right, rotate=0] at (dotrho) {$\dot\varrho$};
			\coordinate (u) at (0.43,0.59) {};
			\node[right, rotate=0] at (u) {$u=\nabla\theta$};
			\coordinate (phi) at (0.84,0.725) {};
			\node[below right, rotate=0] at (phi) {$\varphi$};
			\coordinate (phidot) at (0.92,0.77) {};
			\node[right, rotate=0] at (phidot) {$\dot\varphi$};
			\node[right, rotate=77] at (left_fiber) {fiber};
			\node[right, rotate=77] at (left_diffmu) {$\Diff_\mu$};
			\node[right, rotate=72] at (right_fiber) {fiber};
			\node[above, rotate=5] at (hor_geo) {horizontal geodesic};
			\node[below, rotate=-5] at (geo) {geodesic};
			\node[left, rotate=0] at (projection) {$\pi$};
			\node[above, rotate=-6] at (diff) {$\mathrm{Diff}(M)$};
			\node[above, rotate=-13] at (dens) {$\mathrm{Dens}(M)$};
		\end{scope}
	\end{tikzpicture}
	%
	\caption{Illustration of the Riemannian submersion in \autoref{thm:otto_riemannian_metric}.
	Horizontal geodesics on $\Diff(M)$ (potential solutions) are transversal to the fibres and project to geodesics on $\Dens(M)$.
	}\label{fig:submersion}
\end{figure}
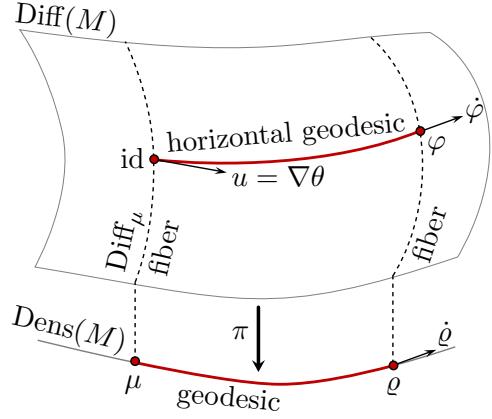

\begin{definition}\label{def:otto_metric}
	The \emph{Wasserstein-Otto metric} on $\Dens(M)$ is
	\begin{equation}\label{eq:otto_metric}
		\MetW_\varrho(\dot\varrho,\dot\varrho) = 
		\int_M \abs{\nabla \theta}^2 \varrho, \qquad \dot\rho = \divv(\rho \nabla \theta),
	\end{equation}
	where $\dot\varrho =\dot\rho \vol \in \Omega_0^n(M)$ is a tangent vector to $\Dens(M)$ at the point $\varrho=\rho\mu$.
The Riemannian distance of this metric is the well-known Wasserstein distance, equal to the minimal $L^2$-cost
of transporting one density to another.
\end{definition}
\begin{theorem}[cf.\ \cite{Ot2001}] \label{thm:otto_riemannian_metric}
The left coset projection $\pi$ is a Riemannian submersion with respect to the $L^2$-metric \eqref{eq:L2met} on $\Diff(M)$ 
and the Wasserstein-Otto metric~\eqref{eq:otto_metric} on $\Dens(M)$.
\end{theorem}
An illustration is given in \autoref{fig:submersion}.
In the Appendix we recall a symplectic interpretation of this construction. 
\begin{theorem}\label{thm:newton_on_dens_hamiltonian_form}
Newton's equation on $\Dens(M)$ for the Wasserstein-Otto metric \eqref{eq:otto_metric} and a potential function $\bar U$
corresponds to Hamilton's equations on $T^\ast\Dens(M)$ 
\begin{equation}\label{eq:ham_eq}
\left\{
  \begin{aligned}
	&\dot\theta + \frac{1}{2}\abs{\nabla\theta}^2 + \frac{\delta \bar U}{\delta\varrho}(\varrho) = 0,\\
	&\dot\varrho + \LieD_{\nabla\theta}\varrho = 0 
\end{aligned} \right. 
\end{equation}
with the Hamiltonian function 
\begin{equation*}
	{H}(\varrho,\theta) = \frac{1}{2}\int_M \abs{\nabla\theta}^2\varrho + \bar U(\varrho). 
\end{equation*}
Solutions to these equations correspond to the potential solutions of Newton's equation \eqref{eq:newton_for_ubar} (or \eqref{eq:reduced_newton_diff}) on $\Diff(M)$.
\end{theorem}
%

\arxivorpnas{
\subsection{Example: Classical mechanics}\label{sub:hamilton_jacobi}
}{
\subsection*{Example: Classical mechanics}
}

Given a $C^\infty$ potential function $V$ on the (finite-dimensional) manifold $M$ we can define 
an associated potential function on the (infinite-dimensional) space of densities $\Dens(M)$ 
\begin{equation}\label{eq:classical_mechanics_potential} 
	\bar U(\varrho) = \int_M V \varrho  .
\end{equation}
\begin{proposition}
	Newton's equation on $\Diff(M)$ (\autoref{thm:newton_for_Ubar}) for a potential of the form in \eqref{eq:classical_mechanics_potential} describes the flow of Newton's equation on $M$ with potential function $V$.
	In particular:
	\begin{itemize}
		\item If $t\mapsto \varphi(t,\cdot)$ is a solution to \eqref{eq:newton_for_ubar}, then for each fixed $x\in M$ the curve $t\mapsto \varphi(t,x)$ satisfies Newton's equation on $M$ with potential $V$.
		\item The vector field $u=\dot\varphi\circ\varphi$ in \eqref{eq:reduced_newton_diff} 
		satisfies the inviscid potential Burgers equation 
		\begin{equation*}
			\dot u + \nabla_u u + \nabla V = 0.
		\end{equation*}
	\end{itemize}
\end{proposition}

\begin{corollary}[cf.\ \cite{KhLe2009}]  \label{cor:ham_jacobi}
The momentum variable $\theta$ in Hamilton's equation \eqref{eq:ham_eq} on $T^*\Dens(M)$ 
satisfies the Hamilton-Jacobi equation for the classical mechanics Hamiltonian on $T^*M$
		$$H(x,p)= \frac{1}{2}\met_x(p^\sharp,p^\sharp) + V(x),$$ 
		where $\sharp$ is the ``musical'' isomorphism defined by the metric $\met$ on $M$.
\end{corollary}
%
%

%

\arxivorpnas{
\subsection{Example: Barotropic fluids}\label{sub:compressible_euler_equations} 
}{
\subsection*{Example: Barotropic fluids}
}

Using \autoref{thm:newton_for_Ubar} one may present the equations of compressible fluids in $M$ 
in an extended Arnold framework with quadratic kinetic energy and without a semidirect group structure.
To this end, consider a potential function of the form 
\begin{equation}\label{eq:barotropic_func} 
	\bar U(\varrho) = \int_M e(\rho)\varrho, 
\end{equation} 
where $e\colon \RR_{>0}\to \RR$ is a function describing the internal energy. 
%
\begin{proposition}[see\  \cite{Eb1975} and  \cite{Sm1979}]
	Newton's equation on $\Diff(M)$ (\autoref{thm:newton_for_Ubar}) for a potential of the form \eqref{eq:barotropic_func} describes a compressible barotropic fluid with internal energy $e$.
	In particular, the form \eqref{eq:reduced_newton_diff} of Newton's equation is given by the compressible Euler equations
	\begin{equation}\label{eq:barotropic} 
	\left\{
		\begin{aligned}
			&\dot u + \nabla_u u + \frac{1}{\rho}\nabla P(\rho) = 0  \\
			&\dot\rho + \divv(\rho u) = 0\,,
		\end{aligned} \right. 
	\end{equation} 
	where the pressure function $P\colon \RR_{> 0}\to \RR$ is given by $P(\rho) = e'(\rho)\rho^2$.
\end{proposition}

\autoref{thm:potential_compress} arises as a combination of this proposition and \autoref{thm:otto_riemannian_metric}.

\begin{remark}
If $M$ is the 2-sphere and $e(\rho)=\rho/2$ we recover the shallow water equations describing the surface motion of an ideal incompressible fluid when the wavelength is large compared to the depth (as in the case of tidal waves).
In this interpretation, $u$ is the surface (horizontal) velocity and $\rho$ is the height of the water.
By using the Nash--Moser inverse function theorem, Hamilton \cite[\S\,III.2.2]{Ha1982} gave a local existence result 
for the shallow water equations.
The above setting allows one to adapt Hamilton's arguments and prove a local existence result 
for any barotropic compressible fluid in the $C^\infty$-setting of tame Fr\'echet spaces. 
This illustrates that our geometric framework can be useful in obtaining analytical results.
\end{remark}

\arxivorpnas{
\section{Fisher-Rao geometry of the space of densities} \label{sec:fisher_rao}
}{
\section*{Fisher-Rao geometry of the space of densities}
}

In this section we introduce a different Riemannian structure on $\Diff(M)$. 
It is given by a Sobolev-class inner product on vector fields and induces on $\Dens(M)$ an infinite-dimensional analogue of the Fisher-Rao information metric.
The setting resembles the relation between the $L^2$-metric \eqref{eq:L2met} and 
Wasserstein-Otto metric \eqref{eq:otto_metric} in the previous section with some notable differences.

As before, let $(M,\met)$ be a compact Riemannian manifold with volume form $\mu$.
Define an $H^1$-metric on $\Diff(M)$ by 
\begin{equation}\label{eq:H1met}
		\Met_{\varphi}(u\circ\varphi,v\circ\varphi) = \frac{1}{4}\int_M \met(-\Delta u,v) \vol + 
		F(u,v),
\end{equation}
where $\Delta$ denotes the Hodge Laplacian on vector fields (cf.~\cite{Mo2015}) and $F(u,v)$ is a positive-definite quadratic form depending only on the vertical (divergence free) components of $u$ and $v$.
\begin{remark}
In the applications below we focus on Newton's equations on $\Dens(M)$ (corresponding to horizontal geodesics on $\Diff(M)$), for which only the first term $\int \met(-\Delta u,v) \vol$ in \eqref{eq:H1met} is relevant.%
\footnote{In \cite{Mo2015} the term $F(u,v)$ is chosen as 
$\sum_{i} \int_M \met(u,e_i)\vol \, \int_M \met(v,e_i)\vol$ where $\{e_i\}$ is any $L^2$-orthogonal basis for 
the space of harmonic vector fields.} 
\end{remark}

\begin{definition}\label{def:fisher_rao_metric}
	The Fisher-Rao \arxivorpnas{(information)}{} metric on $\Dens(M)$ is
	\begin{equation}\label{eq:fisher_rao_metric}
\MetF_\varrho(\dot\varrho,\dot\varrho) = \frac{1}{4}\int_M \left(\frac{\dot\varrho}{\varrho}\right)^2\varrho,
	\end{equation}
	where $\dot\varrho \in \Omega_0^n(M)$ is a tangent vector at the point $\varrho$.
\end{definition}
Consider now the (right coset) projection $\Pi\colon\Diff(M)\to \Dens(M)$ between diffeomorphisms and smooth probability densities given by pull-back of the Riemannian volume form $\Pi(\varphi) = \varphi^*\vol$. 
In analogy with the Riemannian submersion in \autoref{thm:otto_riemannian_metric} we have 
\begin{theorem}\label{thm:FR_riemannian_metric}
The right coset projection $\Pi$ is a Riemannian submersion with respect to 
the $H^1$-metric \eqref{eq:H1met} on $\Diff(M)$ and the Fisher-Rao metric \eqref{eq:fisher_rao_metric} on $\Dens(M)$. 
Furthermore, equipped with this metric $\Dens(M)$ is isometric to a subset of the unit sphere in a Hilbert space 
with a round metric. 
\end{theorem}
\begin{remark} \label{rem:FRiso} 
The isometry between $\Dens(M)$ and a subset of the unit Hilbert sphere described in the last statement of 
\autoref{thm:FR_riemannian_metric} is given by the square root map 
$\varrho \mapsto \sqrt{\rho}$ where $\varrho = \rho \mu$, see \cite{KhLeMiPr2013} for details.
\end{remark} 
%
We point out that the setting of \autoref{thm:FR_riemannian_metric} is quite different from \autoref{thm:otto_riemannian_metric} in one respect.
Namely, the Riemannian metric on $\Diff(M)$ in \autoref{thm:otto_riemannian_metric} is right-invariant 
with respect to $\Diffvol(M)$ and thus automatically descends to the right quotient $\Diff(M)/\Diffvol(M)$.
On the other hand, in \autoref{thm:FR_riemannian_metric} the metric is also right-invariant (under a certain condition  on $F(u,v)$ in \eqref{eq:H1met}),  but nevertheless 
descends to the left quotient $\Diffvol(M)\backslash \Diff(M)$. 
Since the right-invariance property is retained after taking the quotient, the Fisher-Rao metric on $\Dens(M)$ 
remains right-invariant with respect to the action of $\Diff(M)$. 
From this perspective, Fisher-Rao offers a richer geometric structure than Wasserstein-Otto. 
%

%
%
%
%
\begin{theorem}\label{thm:newton_on_dens_fisher}
				%
				%
				%
		%
		Newton's equation on $\Dens(M)$ for the Fisher-Rao metric \eqref{eq:fisher_rao_metric} and a $C^1$ potential function $\bar U\colon \Dens(M)\to \RR$ is
		\begin{equation}\label{eq:EL_eq_dens_fisher}
			\ddot\rho - \frac{\dot\rho^2}{2\rho}
			+ 
			\frac{\delta \bar U}{\delta \varrho}\rho
			= \lambda\rho,
		\end{equation} 
		where $\lambda$ is a Lagrange multiplier for the affine constraint $\int_M\varrho = 1$.
		Solutions correspond to the horizontal solutions of Newton's equation on $\Diff(M)$ for the $H^1$-metric \eqref{eq:H1met} and the potential function \eqref{eq:Ubar}.
\end{theorem}

One can also express \eqref{eq:EL_eq_dens_fisher} as Hamilton's equations on $T^\ast\Dens(M)$, with momentum variable $\theta = \dot\rho/\rho$%
.
\arxivorpnas{
\subsection{Example: \texorpdfstring{$\mu$}{mu}-Camassa-Holm equation} 
}{
\subsection*{Example: \texorpdfstring{$\mu$}{mu}-Camassa-Holm equation} 	
}

The one-dimensional periodic $\mu$CH (also known as $\mu$HS) equation 
\begin{equation} \label{eq:muHS} 
\mu(u_t) - u_{xxt} - 2u_x u_{xx} - uu_{xxx} + 2 \mu(u) u_x = 0, 
\quad 
\mu(u) = \int_{S^1} u \, dx 
\end{equation} 
is a nonlinear PDE closely related to the Camassa-Holm and the Hunter-Saxton equations. 
It describes a model for a director field in the presence of an external (e.g., magnetic) force and 
was derived in \cite{KhLeMi2008} as an Euler-Arnold equation on the group $\Diff(S^1)$ of orientation-preserving circle
diffeomorphisms equipped with a right-invariant $H^1$ Sobolev metric.
%
%
The $\mu$CH equation is known to be bihamiltonian and possess smooth as well as cusped soliton solutions. 
\begin{proposition}[cf.\ \cite{KhLeMi2008, Mo2015}] 
	The $\mu$CH equation \eqref{eq:muHS} is Newton's equation on $\Diff(S^1)$ for the $H^1$-metric \eqref{eq:H1met} and vanishing potential $\bar U\equiv 0$.
	The horizontal (mean zero) solutions $u$ of the $\mu$CH equation describe geodesics of the Fisher-Rao metric~\eqref{eq:fisher_rao_metric} on $\Dens(S^1)$.
\end{proposition} 
\begin{remark}
	The geodesic equation for the $H^1$-metric \eqref{eq:H1met} is sometimes called the EPDiff equation. 
	It is a generalization of $\mu$CH to arbitrary $M$.
	Furthermore, the $H^1$-metric \eqref{eq:H1met} induces a factorization of diffeomorphisms that solves 
	an optimal information-transport problem \cite{Mo2015}, analogous to (but different from) Brenier's factorization \cite{Br1991} in the optimal mass-transport problem.
\end{remark}

\arxivorpnas{
\subsection{Example: infinite-dimensional Neumann problem}\label{sect:Neumann} 
}{
\subsection*{Example: infinite-dimensional Neumann problem}
}

The classical Neumann problem describes the motion of a particle on a sphere $S^n$ in the presence of 
a quadratic potential. 
It is known to be a completely integrable system and to be equivalent (up to time reparametrization) 
to the geodesic motion on an ellipsoid, cf. e.g., \cite{Mo1983}. 

As an infinite-dimensional analogue of the Neumann problem, 
consider Newton's equation on the unit sphere $S^\infty(M) = \{ f \in C^\infty(M)\mid \int_M f^2\vol = 1 \}$ 
with the metric induced from $L^2(M,\mu)$ and with the quadratic potential function
\begin{equation}\label{eq:inf_neumann_potential} 
V(f) 
= 
\frac{1}{2}\int_M \lvert\nabla f\rvert^2 \vol. 
\end{equation}
%
%

\begin{proposition}
	Newton's equation for the infinite-dimensional Neumann problem on $S^\infty(M)$ with potential function \eqref{eq:inf_neumann_potential} is
	\begin{equation}\label{eq:inf_neumann_eq}
		\ddot f -\Delta f = -\lambda f, \qquad \lambda = \int_M (\dot f^2 + f \Delta f)\vol,
	\end{equation} 
	where $\lambda$ is the Lagrange multiplier for the constraint $\int_M f^2\vol = 1$.
	%
\end{proposition}
It turns out that this problem also admits a natural interpretation as a Fisher-Rao Newton's equation on $\Dens(M)$. 
In order to describe it, consider Fisher's information functional on $\Dens(M)$ 
\begin{equation}\label{eq:Fisher_info_func}
I(\varrho) = \frac{1}{8}\int_M \frac{\abs{\nabla \rho}^{2}}{\rho}\vol, 
\qquad \text{where} \quad \rho = \frac{\varrho}{\vol}.
\end{equation}	
%
%
\begin{proposition}\label{pro:newton_fisher_information_neumann}
	The infinite-dimensional Neumann problem \eqref{eq:inf_neumann_eq} corresponds to
	Newton's equation \eqref{eq:EL_eq_dens_fisher} on $\Dens(M)$ with respect to the Fisher-Rao 
	metric and Fisher's information functional \eqref{eq:Fisher_info_func} as a potential function.
	The map $\rho \mapsto \sqrt{\rho} \eqqcolon f$ is a local diffeomorphism between the two representations. 
\end{proposition} 

\begin{remark}
Stationary solutions to the Neumann problem on $S^\infty(M)$ correspond to the principal axes of 
the ellipsoid $\pair{f,-\Delta f}_{L^2} = 1$ and have the natural interpretation as the Laplace eigenfunctions on $M$.  
If $M=\mathbb{T}^4$ is the 4-torus equipped with the (pseudo)-Riemannian Minkowski metric then the stationary solutions 
of the corresponding Minkowski-Neumann problem are solutions of the periodic Klein-Gordon equation 
\begin{equation}\label{eq:klein_gordon}
	\ddot f - \Delta f  = -m^2 f , \qquad m\in\RR.
\end{equation}
This equation describes spinless scalar particles (such as the Higgs boson) and plays a fundamental role 
in quantum field theory. 
The parameter $m$ is interpreted as the particle mass, but in our geometric context it is a Lagrange multiplier 
for the constraint $\int_{\mathbb{T}^4} f^2\vol = 1$.
Through \autoref{pro:newton_fisher_information_neumann} we thereby obtain an interpretation of the Klein-Gordon equation 
as describing stationary potential solutions of a hydrodynamical EPDiff system on $\Diff(\mathbb{T}^4)$. 
This observation may be of some interest in quantum physics.
\end{remark}
%

\arxivorpnas{
\section{Geometric properties of the Madelung transform}\label{sec:madelung}
}{
\section*{Geometric properties of the Madelung transform}	
}
In 1927 
\citet{Ma1927} gave a hydrodynamical formulation of the Schr\"odinger equation. 
Using the setting of the previous sections we can now exhibit a number of surprising geometric properties 
of an important transformation that he introduced. 
\begin{definition}\label{def:madelung} 
Let $\rho$ and $\theta$ be real-valued functions on $M$ with $\rho >0$. 
The \emph{Madelung transform} is the mapping $\Phi\colon(\rho,\theta) \mapsto \psi$ defined by 
\begin{equation}\label{eq:madelung_def} 
	\Phi(\rho,\theta) \coloneqq \sqrt{\rho \ee^{\ii\theta}} . 
\end{equation} 
\end{definition} 
Observe that $\Phi$ is a complex extension of the square root map described 
in \autoref{thm:FR_riemannian_metric} 
(see also \autoref{rem:FRiso}).

\arxivorpnas{
\subsection{Madelung transform as a symplectomorphism} \label{sub:madelung_symplec}
}{
\subsection*{Madelung transform as a symplectomorphism}	
}

As the first property we show that the Madelung transformation induces a symplectomorphism from 
the cotangent bundle of probability densities to the projective space of non-vanishing complex functions.

Let $PC^\infty(M,\CC)$ denote the complex projective space of smooth complex-valued functions on $M$.
We represent its elements as cosets $[\psi]$ of the complex $L^2$-sphere of smooth functions, 
where $\psi'\in[\psi]$ if and only if $\psi' = \ee^{\ii \alpha}\psi$ for some $\alpha\in\RR$.
The space $PC^\infty(M,\CC\backslash \{0\})$ is a submanifold of $PC^\infty(M,\CC)$.
\begin{theorem}\label{thm:madelung_symplectomorphism}
	The Madelung transform \eqref{eq:madelung_def} induces a map 
	\begin{equation}\label{eq:madelung_symplectic}
		\Phi\colon T^*\Dens(M)\to PC^\infty(M,\CC\backslash \{0\}) 
	\end{equation} 
	which is a symplectomorphism (in the Fr\'echet topology of smooth functions) with respect to 
	the canonical symplectic structure of $T^*\Dens(M)$ and the complex projective structure of $PC^\infty(M,\CC)$.
\end{theorem}
The Madelung transform has already been shown to be a symplectic submersion 
from $ T^*\Dens(M)$ to the unit sphere of non-vanishing wave functions by \citet{Re2012}. 
The stronger (symplectomorphism) property stated in \autoref{thm:madelung_symplectomorphism} 
is obtained by considering projectivization $PC^\infty(M,\CC\backslash \{0\})$. 

\medskip

\arxivorpnas{
	\subsection{Example: linear and nonlinear Schr\"odinger equations} 
}{
	\subsection*{Example: linear and nonlinear Schr\"odinger equations}
}

Let $\psi$ be a wavefunction and consider the family of Schr\"odinger (or Gross-Pitaevsky) equations 
(with Planck's constant $\hbar=1$ and mass $m=1/2$) of the form 
\begin{equation}\label{eq:schrodinger} 
 	\mathrm{i}\dot\psi = - \Delta\psi +  V\psi + f(\abs{\psi}^2)\psi, 
\end{equation} 
where $V\colon M\to \RR$ and $f\colon \RR_{> 0}\to \RR$. 
If $f\equiv 0$ we obtain the linear Schr\"odinger equation with potential~$V$. 
If $V\equiv 0$ then we obtain the family of non-linear Schr\"odinger equations (NLS); a typical choice is $f(a) = \kappa a$, another model example is $f(a) = \frac 12(a-1)^2$.

The Schr\"odinger equation \eqref{eq:schrodinger} is a Hamiltonian equation with respect to 
the symplectic structure induced by the complex structure of $L^2(M,\CC)$.
Indeed, 
\arxivorpnas{
if $\inner{\cdot,\cdot}_{L^2(M,\CC)}$ denotes the Hermitian inner product then the real part 
$\pair{\psi,\psi'}_{L^2(M,\CC)} = \mathrm{Re}\inner{\psi,\psi'}_{L^2(M,\CC)}$ 
defines a Riemannian structure and the imaginary part 
\begin{equation*}
	\Omega(\psi,\psi'): = \mathrm{Im}\inner{\psi,\psi'}_{L^2(M,\CC)} = \pair{\mathrm{i}\psi,\psi'}_{L^2(M,\CC)} 
\end{equation*}
defines a symplectic structure.
}{
	the real part of the Hermitian inner product defines a Riemannian structure and the imaginary part defines a symplectic structure.
}
This symplectic form corresponds to the complex structure $J:\psi \mapsto \mathrm{i}\psi$ 
and the Hamiltonian associated with \eqref{eq:schrodinger} is 
\begin{equation*} 
H(\psi) 
= 
\frac{1}{2}\norm{\nabla\psi}_{L^2(M,\CC)}^{2} + \frac{1}{2}\int_M \left( V \abs{\psi}^2 + F(\abs{\psi}^{2}) \right)\vol,
\end{equation*} 
where $F\colon \RR_{> 0}\to \RR$ is a primitive of $f$. 


Observe that the $L^2$ norm of a wave function satisfying \arxivorpnas{the Schr\"odinger equation}{Schr\"odinger's} \eqref{eq:schrodinger} 
is conserved in time. 
Furthermore, the equation is also equivariant with respect to a constant change of phase 
$\psi(x)\mapsto e^{i\alpha}\psi(x)$ and so it descends to the projective space $PC^\infty(M,\CC)$.
Geometrically, the Schr\"odinger equation is thus an equation on complex projective space, 
a point of view first suggested by Kibble \cite{Ki1979}.

%
%
%
%
\begin{proposition}[cf.\ \cite{Ma1927,Re2012}] 
The Madelung transform $\Phi$ 
maps the family of Schr\"odinger equations \eqref{eq:schrodinger} to a family of Newton's equations \eqref{eq:ham_eq} 
on $\Dens(M)$ 
equipped with the Wasserstein-Otto metric \eqref{eq:otto_metric} 
and with potential functions 
\begin{equation}
\bar U(\varrho) = 4 I(\varrho) + \int_M V\varrho + \int_M F(\rho)\vol, 
\end{equation}
where $I$ is Fisher's information functional \eqref{eq:Fisher_info_func}. 
Furthermore, the extension \eqref{eq:reduced_newton_diff} to a system on $\Xcal(M)\times\Dens(M)$ is
\begin{equation}\label{eq:barotropic2} 
\left\{ 
  \begin{aligned} 
	&\dot v + \nabla_v v + \nabla\Big(V + f(\rho) - \frac{2\Delta\sqrt{\rho}}{\sqrt{\rho}} \Big) = 0 \\ 
	&\dot\rho +\divv(\rho v) = 0 \,.
\end{aligned} \right. 
\end{equation} 
\end{proposition}

%
%
%

\begin{corollary}
The Hamiltonian system (\ref{eq:ham_eq}) on $T^*\Dens(M)$ for potential solutions of \eqref{eq:barotropic2} is mapped to 
the Schr\"odinger equation (\ref{eq:schrodinger}) by a symplectomorphism. 
\end{corollary}
Conversely, classical hydrodynamic PDEs can be expressed as NLS-type equations. 
In particular, potential solutions of the compressible Euler equations of a barotropic fluid \eqref{eq:barotropic} 
can be formulated as an NLS equation with Hamiltonian
\begin{equation}\label{eq:compressible_Euler_NLS_Hamiltonian}
H(\psi) 
= 
\frac{1}{2}\norm{\nabla\psi}_{L^2}^{2} 
- 
\frac{1}{2}\norm{\nabla \abs{\psi}}_{L^2}^{2} 
+ 
\int_M e(\abs{\psi}^{2})\abs{\psi}^{2}\vol. 
\end{equation} 
The choice $e=0$ gives a Schrödinger formulation for potential solutions of Burgers' equation (or Hamilton-Jacobi equation, cf.\ \autoref{cor:ham_jacobi}), whose solutions describe geodesics of the Wasserstein-Otto metric \eqref{eq:otto_metric} on $\Dens(M)$. 
Thus the geometric framework connects optimal transport for cost functions with potentials,
the compressible Euler equations and the NLS-type equations described above. 

\medskip


\begin{figure}
	\centering
	\begin{tikzpicture}
		\node[anchor=south west, inner sep=0] (image) at (0,0) {\includegraphics[width=0.25\textwidth]{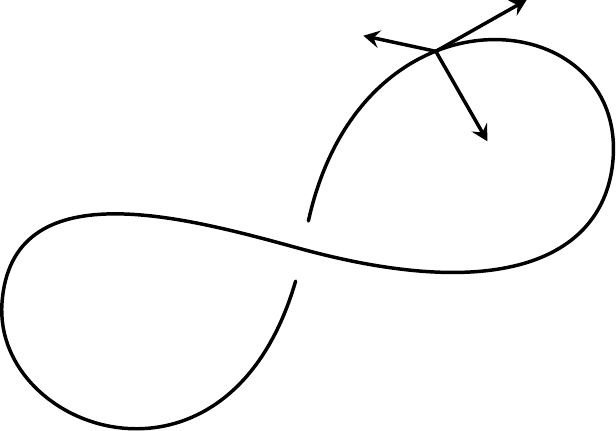}};
		\begin{scope}[x={(image.south east)},y={(image.north west)}]
			\coordinate (gdot) at (0.6,0.92) {};
			\coordinate (gprime) at (0.83,0.93) {};
			\coordinate (gbis) at (0.75,0.73) {};
			\coordinate (g) at (0.45,0.2) {};
			\node[left, rotate=0] at (gdot) {$\dot\gamma = \gamma'\times\gamma''$};
			\node[above right, rotate=0] at (gprime) {$\gamma'$};
			\node[below right, rotate=0] at (gbis) {$\gamma''$};
			\node[right, rotate=0] at (g) {$\gamma$};
		\end{scope}
	\end{tikzpicture}
	\caption{Vortex filament flow.}\label{fig:binormal}
\end{figure}
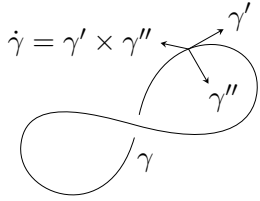

\arxivorpnas{
\subsection{Example: Vortex filament equation} 
}{
\subsection*{Example: Vortex filament equation} 
}
The celebrated vortex filament (or binormal) equation
$$
\dot \gamma=\gamma'\times \gamma''\,,
$$
is an evolution equation for a (closed) curve $\gamma\subset \RR^3$, 
where  $\gamma=\gamma(t,x)$ and $\gamma':=\partial \gamma/\partial x$ and where $x$ is the arc-length parameter 
(see \autoref{fig:binormal}).
It describes a localized induction approximation (LIA) of the 3D Euler equation, where vorticity of the initial velocity field is supported on a curve $\gamma$. 

This  equation is known to be Hamiltonian with respect to the Marsden-Weinstein symplectic structure 
on the space of curves in $\RR^3$ with Hamiltonian given by the length functional, see e.g.\ \cite{ArKh1998}. 
On the other hand, it becomes the equation of the 1D barotropic-type fluid \eqref{eq:barotropic} 
with $\rho=k^2$ and $u=\tau$, where $k$ and $\tau$ denote curvature and torsion of the curve $\gamma$, respectively. 
\smallskip

In 1972 Hasimoto~\cite{Ha1972} discovered the following surprising transformation.
\begin{definition} \label{def:hasimoto} 
The {\it Hasimoto transformation} assigns 
a wave function $\psi:\RR\to\CC$ to a curve $\gamma$ with curvature $k$ and torsion $\tau$, according to the formula 
$$ 
(k(x),\tau(x))\mapsto \psi(x)=k(x)e^{\ii\int^x\tau(\tilde x)d\tilde x}. 
$$ 
\end{definition} 
This map takes the vortex filament  equation to the 1D NLS equation
$\ii\dot \psi+\psi''+\frac 12 |\psi|^2\psi=0\,.$ 
In particular, the filament equation becomes a completely integrable system 
(since the 1D NLS is one) 
whose first integrals are obtained by pulling back those of the NLS equation. 


%
\begin{proposition} The Hasimoto transformation is the 
 Madelung transform in the 1D case.
 \end{proposition}
 This can be seen by comparing Definitions \ref{def:madelung} 
 and \ref{def:hasimoto}, which make the Hasimoto transform seem much less surprising.
 Alternatively, one can note that for $\psi(x)=\sqrt{\rho(x)}e^{\ii\theta(x)/2}$
 the pair $(\rho, u)$ with $u=\nabla\theta$ satisfies the compressible Euler equation, while
 in the one dimensional case these variables are expressed via the curvature $\sqrt{\rho}=\sqrt{k^2}=k$
 and the (indefinite) integral of torsion 
 $\theta(x)/2=\int^x u(\tilde x)d\tilde x=\int^x \tau(\tilde x)d\tilde x$.

\begin{remark} 
The filament equation  has a higher dimensional analog for membranes 
(i.e., compact oriented surfaces $\Sigma$ of codimension 2 in $\RR^n$) 
as a skew-mean-curvature flow
$$
\partial_t q={\bf J}({\bf MC}(q))
$$
where $q\in \Sigma$ is any point of the membrane, ${\bf MC}(q)$ is the mean curvature vector to $\Sigma$ at the point $q$ 
and $\bf J$ is the operator of rotation by $\pi/2$ in the positive direction in every normal space to $\Sigma$. 
This equation is again Hamiltonian with respect to the Marsden-Weinstein structure on membranes of codimension 2 
and with a Hamiltonian function given by the $(n-2)$-dimensional volume of the membrane, see e.g.\ \cite{Sh2012}.

{\bf Open question.} Find an analog of the Hasimoto map, which sends a 
skew-mean-curvature flow to an NLS-type equation for any $n$.
The existence of the Madelung transform and its symplectic property in any dimension 
is a strong indication that such an analog should exist. 
\end{remark}
\medskip

\arxivorpnas{
\subsection{Madelung transform as 
a K\"ahler morphism}\label{sub:kahler_properties_of_madelung}
}{
\subsection*{Madelung transform as 
a K\"ahler morphism}
}

In this section we consider again the Madelung transform as a map between $T^{*}\Dens(M)$ and $P C^{\infty}(M,\CC)$ 
but now equipped with suitable Riemannian structures. 

Consider the tangent bundle $TT^\ast\Dens(M)$ of $T^\ast\Dens(M)$. 
Its elements can be described as 4-tuples 
$(\varrho,\theta,\dot\varrho,\dot\theta)$ 
with $\varrho\in\Dens(M)$, $[\theta] \in C^{\infty}(M)/\RR$, $\dot\varrho \in \Omega^{n}_0(M)$ 
and $\dot\theta \in C^{\infty}(M)$ subject to the constraint 
\begin{equation}
\int_M \dot\theta \varrho = 0. 
\end{equation} 
\begin{definition} 
The 
\emph{Sasaki} (or 
\emph{Sasaki-Fisher-Rao{\rm )} metric} on $T^\ast\Dens(M)$
is the lift of the Fisher-Rao metric \eqref{eq:fisher_rao_metric} from $\Dens(M)$: 
\begin{equation}\label{eq:sasaki_FR_metric}
	\MetF^*_{(\varrho,[\theta])}\left((\dot\varrho,\dot\theta),(\dot\varrho,\dot\theta)\right) = 
	\frac{1}{4}\int_{M} \left( \left(\frac{\dot\varrho}{\varrho}\right)^2 + \dot\theta^2 \right) \varrho .
\end{equation}
The canonical metric on $PC^\infty(M,\CC)$  
\begin{equation} \label{eq:fubini_study}
\Met^*_\psi(\dot\psi,\dot\psi) 
=  
\frac{\pair{\dot\psi,\dot\psi}}{\pair{\psi,\psi}} - \frac{\pair{\psi,\dot\psi}\pair{\dot\psi,\psi}}{\pair{\psi,\psi}^{2}} 
\end{equation} 
is the (infinite-dimensional) \textit{Fubini-Study metric}. 
\end{definition} 
%
%
\begin{theorem}\label{thm:madelung_isometry} 
The Madelung transform~\eqref{eq:madelung_symplectic} is an isometry between $T^\ast\Dens(M)$ equipped with 
the Sasaki-Fisher-Rao metric~\eqref{eq:sasaki_FR_metric} and $P C^{\infty}(M,\CC\backslash \{0\})$ equipped with 
the Fubini-Study metric~\eqref{eq:fubini_study}. 
\end{theorem}
Since the Fubini-Study metric together with the complex structure of $PC^\infty(M,\CC)$ defines a Kähler structure 
it follows that $T^*\Dens(M)$ also admits a natural Kähler structure compatible with its canonical symplectic structure.
\begin{remark}
\citet{Mo2015b} showed that there is an almost complex structure on $T^*\Dens(M)$ 
corresponding to the Wasserstein-Otto metric and the Madelung transform, which does not integrate 
to a complex structure. 
In contrast, our result here shows that the corresponding complex structure becomes integrable (and simple) 
when the Fisher-Rao metric is used instead of Wasserstein-Otto. 
\end{remark} 
%

\arxivorpnas{
\subsection{Example: 2-component Hunter-Saxton (2HS) equation}
}{
\subsection*{Example: 2-component Hunter-Saxton (2HS) equation}
}
The 2HS equation is a system of two equations 
\begin{equation}\label{eq:two_HS_eq}
\left\{
\begin{array}{l}
		\dot u_{xx} = -2 u_x u_{xx} - u u_{xxx} + \sigma\sigma_x, \\
		\dot\sigma = - (\sigma u)_x 
\end{array} \right. 
\end{equation} 
where $u(t,x)$ and $\sigma(t,x)$ are time-dependent periodic functions on the line. 
It can be viewed as a high-frequency limit of the two-component Camassa-Holm equation, cf.\ \cite{HaWu2011}. 

It turns out that this system is closely related to the K\"ahler geometry of the Madelung transformation 
and the Sasaki-Fisher-Rao metric~\eqref{eq:sasaki_FR_metric}. 
Consider the semi-direct product $\mathcal{G} = \Diff_0(S^1)\ltimes C^{\infty}(S^1,S^1)$, 
where $\Diff_0(S^1)$ is the group of circle diffeomorphisms fixing a prescribed point and 
$C^{\infty}(S^1,S^1)$ stands for $S^1$-valued maps of a circle.
Define a right-invariant Riemannian metric on $\mathcal G$ given at the identity by 
\begin{equation*}
\Met_{(\id,0)}\big( (u,\sigma), (v,\tau) \big) = \frac{1}{4} \int_{S^1} \left( u_x v_x +  \sigma\tau\right) \ud x. 
\end{equation*} 
If $t \to (\varphi(t), \alpha(t))$ is a geodesic in $\mathcal{G}$ then 
$u = \dot\varphi\circ\varphi^{-1}$ and $\sigma = \dot\alpha\circ\varphi^{-1}$ 
satisfy equations \eqref{eq:two_HS_eq}, cf.~\cite{Ko2011}.
 Lenells \cite{Le2013b} showed that the map 
\begin{equation}\label{eq:lenells_map} 
(\varphi,\alpha) \mapsto \sqrt{\varphi_x \, \mathrm{e}^{\ii\alpha}} 
\end{equation} 
is an isometry from $\mathcal G$ to an open subset of 
$S^{\infty} = \{ \psi \in C^{\infty}(S^1,\CC) \mid \norm{\psi}_{L^2} = 1 \}. $
Moreover, solutions to \eqref{eq:two_HS_eq} satisfying $\int_{S^1} \sigma \ud x = 0$ correspond to 
geodesics on the complex projective space $PC^{\infty}(S^1,\CC)$ equipped with the Fubini-Study metric.
Our results show that this isometry is a particular case of \autoref{thm:madelung_isometry}. 
%
\begin{proposition}\label{prop:2HS_as_Sasaki}
The 2-component Hunter--Saxton \arxivorpnas{equation}{} \eqref{eq:two_HS_eq} with initial data satisfying 
$\int_{S^1} \sigma \,\ud x = 0$ is equivalent to the geodesic equation of the Sasaki-Fisher-Rao metric 
\eqref{eq:sasaki_FR_metric} on $T^\ast\Dens(S^1)$. 
\end{proposition}
The proof is based on the observation that the mapping \eqref{eq:lenells_map} can be given as 
$(\varphi,\alpha) \mapsto \Phi(\pi(\varphi),\alpha),$
where $\Phi$ is the Madelung transform and $\pi$ is the projection $\varphi\mapsto \varphi^*\vol$ specialized to the case $M=S^1$. 
\begin{remark}
Observe that if $\sigma=0$ at $t=0$ then $\sigma(t)=0$ for all $t$ and 
the 2-component Hunter-Saxton \arxivorpnas{equation}{} \eqref{eq:two_HS_eq} reduces to the standard Hunter-Saxton equation. 
Geometrically, this is a consequence of the fact that horizontal geodesics on $T^*\Dens(M)$ with respect to 
the Sasaki-Fisher-Rao metric descend to geodesics on $\Dens(M)$ with respect to the Fisher-Rao metric.
\end{remark} 
%


\arxivorpnas{
\subsection{Appendix: Madelung transform as a momentum map} \label{sub:madelung_momentum}
}{
\subsection*{Appendix: Madelung transform as a momentum map}	
}
The Riemannian submersion result in \autoref{thm:otto_riemannian_metric} above can be regarded as a Hamiltonian
reduction of the natural symplectic structure on $T^*\Diff(M)$ with respect to the cotangent lifted action of the group $\Diffvol(M)$ 
%
\begin{equation*} 
T^*\Diff(M)\sslash\Diffvol(M) \coloneqq J^{-1}([0])/\Diffvol(M)=T^*\Dens(M) 
\end{equation*}
where $J$ is the associated momentum map, cf.\ \cite{KhLe2009}.

Furthermore, in \autoref{thm:madelung_symplectomorphism} above we described the Madelung transform as a symplectomorphism from $T^*\Dens(M)$ to $PL^2(M,\CC)$.
Following \cite{Fu2017}, we outline here another approach, which shows that the inverse Madelung map 
is a momentum map from the space  $\Psi\subset C^\infty(M,\CC)$ of smooth 
complex-valued wave functions $\psi$ with the $L^2$-product to the set of pairs $(\rho, \nabla\theta)$,  
regarded as elements of the dual space $\mathfrak{s}^*$ to the semidirect product 
Lie algebra $\mathfrak s=\Xcal(M)\ltimes C^\infty(M)$ of the group $S=\Diff(M)\ltimes C^\infty(M)$. 

It is convenient to think of $\Psi$ as a space of complex-valued \emph{half-densities} on $M$.
Half-densities are characterized by how they are transformed under diffeomorphisms of the underlying space: 
the pushforward $\varphi_*\psi$ of a half-density $\psi$ on $M$ by a diffeomorphism $\varphi\in\Diff(M)$ is given by
$$
\varphi_* \psi = \sqrt{|\mathrm{Det}(D\varphi^{-1})|}\, \psi \circ \varphi^{-1}\,.
$$
This explains the following natural action of $S$ on the space $\Psi$ interpreted as half-densities.
\begin{definition}[\cite{Fu2017}] 
The semidirect product 
	$S=\mathrm{Diff}(M)\ltimes C^\infty(M)$ acts on
	$\Psi$ by
\begin{equation*}
		(\varphi,a) \cdot \psi = \sqrt{|\mathrm{Det}(D\varphi^{-1})|}\, e^{-\ii a} (\psi \circ \varphi^{-1}).
\end{equation*}
\end{definition}
While the classical Madelung transform $(\theta,\rho) \mapsto \sqrt{\rho e^{\ii\theta}}$ is defined for positive $\rho$, 
the inverse map has a particularly nice form.
Assume for now that $M=\RR^n$.

\begin{proposition}
	The map  ${\bf M} : \psi\mapsto(\mu, \rho):= ({\rm Im} \, \bar \psi \nabla \psi ,\bar \psi \psi )$ is the inverse of the classical Madelung transform in the following sense.
	If $\psi = \sqrt{\rho e^{\ii\theta}}$, then ${\bf M}(\psi) = (\rho \nabla \theta , \rho)$. 
	If $\rho>0$, then the pair $(\rho\nabla \theta, \rho)$ can be identified with $(\rho,[\theta])\in T^*\Dens(M)$.
\end{proposition}

This follows from the observation that 
$\text{Im} \, \bar \psi \nabla \psi=
 \bar \psi \psi\, \text{Im} \,  \nabla \ln \psi $.

\begin{theorem}[\cite{Fu2017}]
The action of the group $S$ on the space $\Psi\subset C^\infty(M,\CC)$ preserves the symplectic structure 
on $\Psi$ and is Hamiltonian.
The corresponding momentum map $ \Psi \to \mathfrak s^*$ is the inverse Madelung transform ${\bf M}$.\end{theorem}

In particular, this theorem implies that  the Madelung transform is also a Poisson map taking the bracket 
on $\Psi$ to the Lie-Poisson bracket on $\mathfrak s^*$.




\IfFileExists{/Users/moklas/Documents/Papers/References.bib}{
\bibliographystyle{amsplainnat}
\bibliography{/Users/moklas/Documents/Papers/References} 
}{
\def\cprime{$'$}

}
 
\end{document}